\documentclass[a4paper,12pt]{amsart}

\usepackage{graphicx}
\usepackage{amsmath,amssymb,amsthm,amsfonts}
\usepackage{amssymb}

\newtheorem{thm}{Theorem}[section]

\newtheorem{lem}[thm]{Lemma}

\newtheorem{rem}[thm]{Remark}

\newcommand{\bremark}{\begin{rem} \textup}
\newcommand{\eremark}{\end{rem} }

\newcommand{\cuad}{{\sqcap\kern-.68em\sqcup}}

\newcommand{\Om}{{\Omega}}

\newcommand{\R}{{\mathbb{R}}}

\begin{document}

\parindent 0pc
\parskip 6pt
\overfullrule=0pt

\title[Extremal solution for singular p-Laplace equations]{Regularity of the extremal solution for singular p-Laplace equations}
       \date{July 1, 2014}
\author{Daniele Castorina}
\thanks{Daniele Castorina, Dipartimento di Matematica, Universit\'a di Padova,
Via Trieste 63, 35121 Padova, Italy. e-mail: castorin@math.unipd.it}

\begin{abstract}
We study the regularity of the extremal solution $u^*$ to the singular
reaction-diffusion problem $-\Delta_p u = \lambda f(u)$ in
$\Omega$, $u =0$ on $\partial \Omega$, where $1<p<2$, $0 < \lambda <
\lambda^*$, $\Omega \subset \mathbb{R}^n$ is a smooth bounded domain
and $f$ is any positive, superlinear, increasing and (asymptotically)
convex $C^1$ nonlinearity. We provide a simple proof of known $L^r$ and $W^{1,r}$
\textit{a priori} estimates for $u^*$, i.e.  $u^* \in L^\infty(\Omega)$ if $n \leq p+2$, $u^* \in L^{\frac{2n}{n-p-2}}(\Omega)$ if $n > p+2$ and $|\nabla u^*|^{p-1} \in
L^{\frac{n}{n-(p'+1)}} (\Omega)$ if $n > p p'$. 
\end{abstract}

\maketitle


\section{Introduction and main result}\label{intro}

\setcounter{equation}{0}

The aim of this paper is the study of the following quasilinear
reaction-diffusion problem:

\begin{equation}\label{problem}
\left\{
\begin{array}{rcll}
-\Delta_p u &=& \lambda f(u)  &\textrm{in } \Omega, \\
 u&>& 0  &\textrm{in } \Omega, \\
 u &=& 0  &\textrm{on } \partial \Omega,
\end{array}
\right.
\end{equation}

where the diffusion is driven by the singular p-Laplace operator $\Delta_p
u:= {\rm div}(|\nabla u|^{p-2}\nabla u)$, $1<p<2$, $\Omega$ is a
smooth bounded domain of $\R^n$, $n\geq 2$, $\lambda$ is a positive
parameter and the nonlinearity $f$ is any $C^1$ positive increasing
function satisfying

\begin{equation}\label{p-superlinear}
\lim_{t\rightarrow+\infty}\frac{f(t)}{t^{p-1}}=+\infty.
\end{equation}

Typical reaction terms $f$ satisfying the above assumptions are
given by the exponential $e^u$ and the power $(1+u)^m$ with $m>p-1$.

These reaction-diffusion problems appear in numerous models in
physics, chemistry and biology. In particular, when $p=2$ and $f(u)$
is the exponential, $\eqref{problem}_{\lambda}$ is usually referred
to as the \textit{Gelfand problem}: it arises as a simplified model
in a number of interesting physical contexts. For example, up to
dimension $n=3$, equation $\eqref{problem}_{\lambda}$ can be derived
from the thermal self-ignition model which describes the reaction
process in a combustible material during the ignition period. We
refer the interested reader to \cite{BE,FK} for the detailed
derivation of the model, as well as other physical motivations for
this problem. In the case of singular nonlinearities such as
$f(u)=(1-u)^{-2}$, problem $\eqref{problem}_{\lambda}$ is also
relevant as a model equation to describe Micro Electro Magnetic
System (MEMS) devices theory (see \cite{PeBe} for a complete account
on this subject). Regarding the MEMS equation, the compactness of
the minimal branch of solutions and some spectral issues connected
with it were investigated in \cite{CES1} for general $p>1$ and
nonlinearities $f(u)$, singular at $u=1$, with growth comparable to $(1-u)^{-m}$, $m>0$.\\

In order to set up the problem, we will say that a (nonnegative)
function $u\in W^{1,p}_0(\Omega)$ is a \emph{weak energy solution}
of $\eqref{problem}_{\lambda}$ if $f(u)\in L^1(\Omega)$ and $u$
satisfies

$$
\int_\Omega |\nabla u|^{p-2}\nabla u\cdot\nabla \varphi\ dx= \lambda
\int_\Omega f(u)\varphi\ dx\quad \textrm{for all } \varphi\in
C_0^1(\Omega).
$$

Moreover, if $f(u)\in L^\infty(\Omega)$ we say that $u$ is a
\emph{regular solution} of $\eqref{problem}_{\lambda}$. By standard
regularity results for non-uniformly elliptic equations, one has
that every regular solution $u$ belongs to $C^{1,\alpha}
(\overline{\Omega})$ for some $0<\alpha<1$ (see \cite{DB,Lie,T}).\\

Under the above hypotheses, problem $\eqref{problem}_{\lambda}$ has
been extensively studied for $p=2$. Crandall and Rabinowitz
\cite{CR} prove  the existence of an extremal parameter
$\lambda^*\in (0,+\infty)$ such that: if $\lambda<\lambda^*$ then
problem $\eqref{problem}_{\lambda}$ admits a regular solution
$u_\lambda$ which is minimal among all other possible solutions, and
if $\lambda>\lambda^*$ then problem $\eqref{problem}_{\lambda}$
admits no regular solution. Moreover, they prove that every minimal
solution $u_\lambda$ is semi--stable in the sense that the second
variation of the energy functional associated to
$\eqref{problem}_{\lambda}$ is nonnegative definite. Subsequently,
Brezis and V\'azquez \cite{BreVaz97} prove that the pointwise
increasing limit of minimal solutions given by

\begin{equation}\label{eq5}
u^*:=\lim_{\lambda\uparrow\lambda^*}u_\lambda,
\end{equation}

is a weak solution of $\eqref{problem}_{\lambda}$, usually known as
\textit{extremal solution}. Apart from a detailed study of the model cases
(exponential and power nonlinearities), they also raise some
interesting open problems for general nonlinearities $f$ satisfying
the above assumptions. One of the most challenging questions
proposed, which has not been answered completely yet, is to show
show that the extremal solution $u^*$ is bounded or in the energy
class, depending on the range of dimensions. In this direction,
Nedev \cite{Nedev} proved, in the case of convex nonlinearities,
that $u^*\in L^\infty(\Omega)$ if $n\leq 3$ and $u^*\in
L^r(\Omega)$ for all $1\leq r<n/(n-4)$ if $n\geq 4$. Subsequently,
Cabr\'e~\cite{Cabre09}, Cabr\'e and Sanch\'on \cite{CS}, and Nedev
\cite{Nedev01} proved, in the case of convex domains and general
nonlinearities, that $u^*\in L^\infty(\Omega)$ if $n\leq 4$ and
$u^*\in L^{\frac{2n}{n-4}}(\Omega)\cap H^1_0(\Omega)$ if $n\geq
5$. More recently, Villegas \cite{Vi} extends Nedev result to $n=4$
thanks to a clever use of the \textit{a priori} estimates of
Cabr\'e~\cite{Cabre09}, without convexity assumptions on the domain.\\

While for the standard Laplacian the regularity of the extremal
solution has been subject of a rich literature, in the case of the
p-Laplacian, i.e. for $p\ne2$, the available results are very few.
Garc\'{\i}a Azorero, Peral and Puel \cite{GarPe92,GarPePu94} study
$\eqref{problem}_{\lambda}$ when $f(u)=e^u$, obtaining the existence
of the family of minimal regular solutions $(\lambda,u_\lambda)$ for
$\lambda\in (0,\lambda^*)$ and that $u^*$ is a weak energy solution
independently of $n$. If in addition $n<p+4p/(p-1)$ then $u^*\in
L^\infty(\Om)$. Moreover, they prove that $\lambda^*=p^{p-1}(n-p)$
and $u^*(x)={\rm log}(1/|x|^p)$ when $\Omega=B_1$ and $n\geq
p+4p/(p-1)$. Cabr\'e and Sanch\'on \cite{CS07} proved the existence
of an extremal parameter $\lambda^\star\in(0,\infty)$ such that
problem $\eqref{problem}_{\lambda}$ admits a minimal regular
solution $u_\lambda\in C^1_0(\overline{\Omega})$ for
$\lambda\in(0,\lambda^*)$ and admits no regular solution for
$\lambda>\lambda^*$. Moreover, every minimal solution $u_\lambda$ is
semi-stable for $\lambda\in (0,\lambda^*)$.

One of the main difficulties is the fact that for arbitrary $p>1$ it
is unknown if the limit of minimal solutions $u^*$ is a (weak or
entropy) solution of $\eqref{problem}_{\lambda^*}$. In the
affirmative case, it is called the \textit{extremal solution} of
$\eqref{problem}_{\lambda^*}$. However, in \cite{S} Sanch\'on has
proved that the limit of minimal solutions $u^*$ is a weak solution
(in the distributional sense) of $\eqref{problem}_{\lambda^*}$
whenever $p\geq 2$ and $(f(t)-f(0))^{1/(p-1)}$ is convex for $t$
sufficiently large. Essentially under the same hypotheses,
Bidaut-Veron and Hamid in \cite{BH} are able to show that $u^*$ is a
locally renormalized solution of $\eqref{problem}_{\lambda^*}$ in the singular case $p < 2$.\\

In this paper, in the spirit of the clever proof of \cite{Vi} for the case $p=2$, we extend some of the results of
\cite{CaSa} for the degenerate p-Laplacian to the singular case $1<p<2$. We obtain the boundedness of
the extremal solution up to a critical dimension $n_p = p+2$ while we prove that it belongs to $L^{\frac{2n}{n-p-2}}(\Omega)$ if $n > p+2$, for any smooth bounded domain $\Omega$, 
under a standard (asymptotic) convexity assumption on the nonlinearity. Unfortunately our dimensional $n_p$ is not optimal in the singular case ($n_p < p p'$ for $1<p<2$), hence we are not able to match the one obtained in \cite{BH} for this case. However, the regularity results in \cite{BH} for the singular case are rather involved while our alternative proof is very simple and direct. In higher dimensions, we will establish in a different way the same Sobolev regularity obtained in \cite{S} for the degenerate case. Our main result is the following:\\

\begin{thm}\label{teorema}
Suppose that $1<p<2$ and let $f$ be a positive, increasing and
superlinear $C^1$ nonlinearity such that $f^{\frac{1}{p-1}}(t)$ is
convex for any $t \geq T$. Let $u^*$ be the extremal solution of
\eqref{problem} and set $n_p := p+2$. The following assertions hold:

$(a)$ If $n \leq n_p$ then $u^* \in L^\infty(\Omega)$. In particular,
$u^*$ is a regular solution to $\eqref{problem}_{\lambda^*}$.

$(b)$ If $n > n_p$ then $u^* \in L^{\frac{2n}{n-p-2}}(\Omega)$.

$(c)$ If $n> p p'$ then $|\nabla u^*|^{p-1} \in L^{\frac{n}{n-(p'+1)}} (\Omega)$.
\end{thm}

We will recall several auxiliary results as well as giving the proof of Theorem \ref{teorema} in the next section. The main details of the proof of the technical lemmas can be found in the Appendix, section \ref{appendix}.

\section{Proof of Theorem \ref{teorema}}\label{proof}

\setcounter{equation}{0}

Let us discuss a few preliminary results which will be used. First
of all, the proof of Theorem \ref{teorema} relies on the
semistability of the minimal solution $u_\lambda$ for $0 < \lambda <
\lambda^* $.

Recall that the linearization $L_{u_{\lambda}}$ associated to
$\eqref{problem}_\lambda$ at a given solution $u_{\lambda}$ is
defined as

\begin{eqnarray*}
L_{u_{\lambda}} (v,\varphi) &:=& \int_\Omega |\nabla u_\lambda|^{p-2}(\nabla v,\nabla\varphi)\\
&&+(p-2)\int_\Omega |\nabla u_\lambda|^{p-4}(\nabla u_\lambda,\nabla
v)(\nabla u_\lambda,\nabla\varphi) - \int_\Omega \lambda f'(u_\lambda)v
\varphi. \nonumber
\end{eqnarray*}

for test functions $v, \varphi \in C^1_{c} (\Omega)$. Observe that the above linearization, in the degenerate case $1<p<2$, makes sense 
if $|\nabla u_\lambda|^{p-2} \in L^1(\Omega)$, which has been proved by Damascelli and Sciunzi in \cite{DS}.
We then say that a solution of $\eqref{problem}_\lambda$ is {\em semistable} if
the linearized operator at $u_{\lambda}$ is nonnegative definite,
i.e. $L_{u_{\lambda}}(\varphi,\varphi) \geq 0$ for any $\varphi \in
C^1_{c} (\Omega)$. Equivalently, $u_{\lambda}$ is semistable if the
first eigenvalue of $L_{u_{\lambda}}$ in $\Omega$, $\mu_1
(L_{u_{\lambda}},\Omega)$, is nonnegative. However, let us observe
that for $p \ne 2$ the latter definition of semistability requires
the spectral theory for $L_u$ that has been established by Esposito,
Sciunzi and the author in \cite{CES4}.

Next, we will need two \textit{a priori} estimates for the family of
minimal solutions $u_{\lambda}$, $0 < \lambda < \lambda^*$. For the reader's convenience a sketch of the proofs of both auxiliary results can be found in the Appendix (section \ref{appendix}).

The first estimate gives uniform $L^\infty$ and $L^r$ bounds for $u_\lambda$ in
terms of the $W^{1,p+2}_0$ norm in a neighborhood of the boundary, depending on the dimension. 
It can be directly derived from the \textit{a priori}
estimates for semistable solutions contained in Theorem 1.4 of
\cite{CaSa}, where Sanch\'on and the author extend to the case $p>2$
the regularity results of \cite{Cabre09,CS} for $u^*$ in convex
domains.

\begin{lem}\label{lemma1}
Let $u_\lambda$ be the minimal solution of $\eqref{problem}_\lambda$. Then the following alternatives hold:

$(a)$ If $n\leq p+2$ then there exists a constant $C$ depending
only on $n$ and $p$ such that

\begin{equation}\label{L-infinty}
\|u_\lambda\|_{L^\infty(\Omega)}\leq s+\frac{C}{s^{2/p}}|\Omega|^\frac{p+2-n}{np}
\left(\int_{\{u_\lambda\leq s\}}  |\nabla u_\lambda|^{p+2}\, dx\right)^{1/p}\quad \textrm{for
all }s>0.
\end{equation}

$(b)$ If $n>p+2$ then there exists a constant $C$ depending
only on $n$ and $p$ such that

\begin{equation}\label{Lq:estimate}
\left(\int_{\{u_\lambda>s\}} \Big(u_\lambda-s\Big)^{\frac{np}{n-(p+2)}}\ dx\right)^{\frac{n-(p+2)}{np}}
\leq \frac{C}{s^{2/p}}
\left(\int_{\{u_\lambda\leq s\}} |\nabla u_\lambda|^{p+2} \ dx\right)^{1/p}
\end{equation}

for all $s>0$.
\end{lem}

The second preliminary result is a uniform $L^1$ estimate for
$(f(u_\lambda))^{p'}/u_\lambda$ when the power $\frac{1}{p-1}$ of
the nonlinearity is asymptotically convex. This bound is a direct
consequence of the estimates which have been proved in Proposition 5.28
of \cite{BH}.

\begin{lem}\label{lemma2}
Let $u_\lambda$ be the minimal solution of $\eqref{problem}_\lambda$
and suppose that $f^{\frac{1}{p-1}}(t)$ is convex for any
sufficiently large $t$. Then there exists a constant $M$ independent
of $\lambda \in (0,\lambda^*)$ such that

\begin{equation}\label{L-p'}
\int_{\{ u_\lambda > 1 \}} \frac{(f(u_\lambda))^{p'}}{u_\lambda} \,
dx \leq M.
\end{equation}
\end{lem}

Finally, we will need a uniform gradient estimate which we
can derive from the regularity results for the linear problem. Let
us consider

\begin{equation}\label{eq8}
\begin{cases}
-\Delta_p u = g(x)\quad \textrm{ in }\Omega\\
\qquad u=0\qquad \textrm{ on }\partial \Omega
\end{cases}
\end{equation}

where $g\in L^m(\Omega)$ for some $m>1$. The following result gives
the regularity of the gradient of a weak energy solution of
(\ref{eq8}) (see for instance \cite{ABFOT}).

\begin{lem}\label{lem3}
Suppose that $n>q$ and let $q^*:=\frac{nq}{n-q}$ be the usual
critical Sobolev exponent. Let $u$ be a weak energy solution of
\eqref{eq8}. Then there exists a constant $C$ depending on $n$, $p$,
$q$ and $|\Omega|$ such that:$$\||\nabla u|^{p-1}\|_{L^{q^*}
(\Omega)} \leq C \|g\|_{L^q (\Omega)}.$$
\end{lem}

\begin{rem} Notice that for $p<2$ it might happen that $(p-1) q^* <1$. Hence the Sobolev norm
$W_{0}^{1,(p-1)q^*}$ of the solution $u$ might not make sense:
however, for the sake of simplicity, we will still keep this
notation intending throughout the discussion that: $$\| u
\|_{W^{1,(p-1)q^*}_0} := \| |\nabla u|^{p-1}
\|_{L^{q^*}}^{\frac{1}{p-1}}$$
\end{rem}

We are now ready to prove the regularity statements for $u^*$ given
in Theorem \ref{teorema}.\\

{\bf Proof of Theorem \ref{teorema} part (a).}\\

Let us begin by noticing that a trivial consequence of \eqref{L-infinty} is

\begin{equation}\label{L-infinty2}
\|u_\lambda \|_{L^\infty(\Omega)} \leq
s+\frac{C}{s^{2/p}}|\Omega|^\frac{p+2-n}{np} \left(\int_{\Omega} |\nabla u_\lambda |^{p+2}\, dx\right)^{1/p}
\end{equation}

where the constant $C$ is given by Lemma \ref{lemma1}. Setting $A:= \| u \|_{W^{1,p+2}_0 (\Omega)}$, for any
$s>0$ consider the RHS of \eqref{L-infinty2} as given by the function

$$\Phi (s) := s+C A^{\frac{p+2}{p}} s^{-\frac{2}{p}}.$$

By explicit computation we see that

$$\Phi'(s) = 1 - \frac{2C}{p} A^{\frac{p+2}{p}} s^{-\frac{p+2}{p}}$$

Notice that $\Phi$ is a strictly convex function with a unique
global minimum at

$$s = \left(\frac{2C}{p} \right)^{\frac{p}{p+2}} A.$$

By direct substitution we see that

$$\Phi\left(\left(\frac{2C}{p}\right)^{\frac{p}{p+2}} A \right) = \left[\left(\frac{2C}{p}\right)^{\frac{p}{p+2}} + \frac{C^{\frac{4-p}{2-p}}
p^\frac{2}{p-2}}{2^\frac{2}{p-2}}\right] A$$

In particular, by the estimate \eqref{L-infinty2} we have thus deduced that
for any $0 < \lambda < \lambda^*$ there exists a positive constant
$D$ independent of $\lambda$ such that we have

\begin{equation}\label{a1}
\|u_\lambda \|_{L^\infty(\Omega)}\leq D \| u_\lambda \|_{W^{1,p+2}_0
(\Omega)}
\end{equation}

Now let us observe that since $u_\lambda$ is weak energy solution of
$\eqref{problem}_\lambda$, from Lemma \ref{lem3} we know that there
exists a positive constant $E$ such that:

\begin{equation}\label{a2}
\| u_\lambda \|_{W^{1,p+2}_0 (\Omega)} \leq E \| \lambda f(u_\lambda)
\|_{L^{q} (\Omega)}^{\frac{1}{p-1}}
\end{equation}

for $q = \frac{n(p+2)}{(p-1)n+p+2}$. Thus, taking into account
\eqref{L-p'}, \eqref{a1} and \eqref{a2} and recalling that $f$ is
increasing, we obtain that for any $0 < \lambda <
\lambda^*$ the following chain of inequalities holds:

$$
\|u_\lambda \|_{L^\infty(\Omega)}^{q(p-1)} \leq D^{q(p-1)} \|
u_\lambda \|_{W^{1,p+2}_0 (\Omega)}^{q(p-1)} \leq (D E)^{q(p-1)} \|
\lambda f(u_\lambda) \|_{L^{q} (\Omega)}^q
$$

$$
\leq (D E)^{q(p-1)} (\lambda^*)^q \left( \int_{\{u_\lambda \leq 1\}}
(f(u_\lambda))^q \, dx + \int_{\{ u_\lambda > 1 \}} (f(u_\lambda))^q
\, dx \right)
$$

$$
\leq C \left( (f(1))^q |\Omega| + \int_{\{ u_\lambda > 1 \}}
\frac{f(u_\lambda)^q}{u_{\lambda}^{\frac{q}{p'}}}
u_{\lambda}^{\frac{q}{p'}} \, dx \right)
$$

$$
\leq C \left( (f(1))^q |\Omega| + \|u_\lambda
\|_{L^\infty(\Omega)}^{\frac{q}{p'}} |\Omega|^{\frac{q}{p'-q}}
\left(\int_{\{ u_\lambda > 1 \}}
\frac{(f(u_\lambda))^{p'}}{u_\lambda} \, dx\right)^{\frac{q}{p'}}
\right)
$$

$$
\leq C \left( (f(1))^q |\Omega| + M^{\frac{q}{p'}}
|\Omega|^{\frac{q}{p'-q}} \|u_\lambda
\|_{L^\infty(\Omega)}^{\frac{q}{p'}} \right)
$$

with $C:=(D E)^{q(p-1)} (\lambda^*)^q$. Let us point out that in the third line of the above calculations we have applied Holder
inequality under the condition $q < p'$, which is true for $n < q_p := \frac{p(p+2)}{2(p-1)}$ and in fact $n \leq n_p < q_p$ for $1<p<2$. 

Therefore there exist two positive constants $A$ and $B$ independent of $\lambda$ such that:

$$
\|u_\lambda \|_{L^\infty(\Omega)}^{q(p-1)} \leq A + B \|u_\lambda
\|_{L^\infty(\Omega)} ^{\frac{q}{p'}}
$$

Observing that $q(p-1) > q/p'$ for $p>1$, the above inequality
implies that $u_\lambda$ is uniformly bounded in $L^\infty(\Omega)$
and, taking the limit as $\lambda \uparrow \lambda^*$, we see that
$u^* \in L^\infty (\Omega)$. This concludes
the proof of statement (a).\\

{\bf Proof of Theorem \ref{teorema} part (b).}\\

Observe that thanks to \eqref{Lq:estimate} for any $s>0$ we have:

$$\int_{\Omega} |u_\lambda|^{\frac{np}{n-(p+2)}} \, dx =  \int_{\{u_\lambda \leq s\}} |u_\lambda|^{\frac{np}{n-(p+2)}} \, dx + \int_{\{u_\lambda>s\}} |u_\lambda|^{\frac{np}{n-(p+2)}}  \, dx$$
$$= \int_{\{u_\lambda \leq s\}} |u_\lambda|^{\frac{np}{n-(p+2)}} \, dx + \int_{\{u_\lambda>s\}} |(u_\lambda-s)+s|^{\frac{np}{n-(p+2)}} \, dx$$
$$\leq s^{\frac{np}{n-(p+2)}} |\{u_\lambda \leq s\}| + s^{\frac{np}{n-(p+2)}} |\{u_\lambda>s\}| + \int_{\{u_\lambda>s\}} (u_\lambda-s)^{\frac{np}{n-(p+2)}} \, dx$$
$$\leq s^{\frac{np}{n-(p+2)}} |\Omega| + \frac{C}{s^{\frac{2n}{n-(p+2)}}} \left(\int_{\{u_\lambda\leq s\}} |\nabla u_\lambda|^{p+2} \ dx\right)^{\frac{n}{ n-(p+2)}} $$

From the above chain of inequalities we easily deduce that 

\begin{equation}\label{Lq:estimate2}
\|u_\lambda \|_{L^{\frac{np}{n-(p+2)}}(\Omega)}^{\frac{np}{n-(p+2)}} \leq
s^{\frac{np}{n-(p+2)}} |\Omega| +\frac{C}{s^{\frac{2n}{n-(p+2)}}} \left(\int_{\Omega} |\nabla u_\lambda |^{p+2}\, dx\right)^{\frac{n}{ n-(p+2)}}
\end{equation}

At this point, we note that the RHS of the above inequality \eqref{Lq:estimate2}, up to multiplicative constant depending only on $|\Omega|$ and a change of variables $t = s^{\frac{np}{n-(p+2)}}$, is given by the same function $\Phi (t)$ which has been optimized at the beginning of the previous proof. Then, in order to prove part $(b)$ of the theorem, we can proceed essentially as in the proof of part $(a)$ using \eqref{Lq:estimate2} in place of \eqref{L-infinty2}.\\

{\bf Proof of Theorem \ref{teorema} part (c).}\\

Let us now suppose that $n > p p'$. Applying
again Lemma \ref{lem3}, but this time with exponent $q =
\frac{n(p-1)}{n(p-1)-p}$ (which spells $q^*=\frac{n}{n-(p'+1)}$), we see that:

$$
\| u_\lambda \|_{W^{1,(p-1)q^*}_0 (\Omega)} \leq E \| \lambda
f(u_\lambda) \|_{L^{q} (\Omega)}^{\frac{1}{p-1}}
$$
Proceeding exactly as in the proof of part $(a)$ and applying Holder
inequality (notice that for $n > p p'$ and $p < 2$ we have that $q <
p'$) and Sobolev inequality we obtain:
$$
\| u_\lambda \|_{W^{1,(p-1)q^*}_0 (\Omega)}^{q(p-1)} \leq E^{q(p-1)}
(\lambda^*)^q \left( (f(1))^q |\Omega| + \int_{\{ u_\lambda > 1 \}}
\frac{f(u_\lambda)^q}{u_{\lambda}^{\frac{q}{p'}}}
u_{\lambda}^{\frac{q}{p'}} \, dx\right)
$$

$$
\leq E^{q(p-1)} (\lambda^*)^q \left( (f(1))^q |\Omega| +
M^{\frac{q}{p'}} \left(\int_\Omega u_{\lambda}^{\frac{q}{p'-q}}
\right)^{\frac{p'-q}{p'}} \right)
$$

$$
\leq E^{q(p-1)} (\lambda^*)^q \left( (f(1))^q |\Omega| +
M^{\frac{q}{p'}} S^{} \| u_\lambda \|_{W^{1,(p-1)q^*}_0
(\Omega)}^{\frac{q}{p'}} \right)
$$

Once again there exist two positive constants $A$ and $B$
independent of $\lambda$ such that:

$$
\| u_\lambda \|_{W^{1,(p-1)q^*}_0 (\Omega)}^{q(p-1)} \leq A + B \|
u_\lambda \|_{W^{1,(p-1)q^*}_0 (\Omega)}^{\frac{q}{p'}}
$$

We thus see that $u_\lambda$ is uniformly bounded in $W^{1,(p-1)q^*}_0
(\Omega)$ and, taking the limit as $\lambda \uparrow \lambda^*$, we
get that $u^* \in W^{1,(p-1)q^*}_0 (\Omega)$. This is exactly
statement (c) of Theorem \ref{teorema}, so the proof is done. \qed

\section{Appendix}\label{appendix}

{\bf Sketch of the proof of Lemma \ref{lemma1}.}

Let us recall that the semistabilty of $u_\lambda$ reads as 

\begin{equation}\label{semistab}
\int_\Omega |\nabla u_\lambda|^{p-2}|\nabla \varphi|^2 +(p-2)\int_\Omega |\nabla u_\lambda|^{p-4}(\nabla u_\lambda,\nabla\varphi)^2 - \int_\Omega \lambda f'(u_\lambda)\varphi^2 \geq 0
\end{equation}
for any $\varphi \in C^{1}_{c} (\Omega)$. Considering $\phi = |\nabla u_\lambda| \eta$ as a test function in \eqref{semistab}, we obtain
\begin{equation}\label{StZu}
\int_{\Omega}\left[ (p-1) |\nabla u_\lambda|^{p-2} |\nabla_{T,u_\lambda} |\nabla u_\lambda||^{2}
+ B_{u_\lambda}^2 |\nabla u_\lambda|^{p} \right] \eta^2 \, dx
\leq (p-1) \int_{\Omega}  |\nabla u_\lambda|^{p} |\nabla \eta|^2 \, dx
\end{equation}
for any Lipschitz continuous function $\eta$ with compact support. Here $\nabla_{T,v}$ is the tangential gradient along a level set of $|v|$ while $B_v^2$ denotes the $L^2$-norm of the second fundamental form of the level set of $|v|$ through $x$. The fact that $\phi=|\nabla u_\lambda|\eta $ is an admissible test function as well as the computations behind \eqref{StZu} can be found in \cite{FSV0} (see also Theorem 1 in \cite{FSV}).

On the other hand, noting that $(n-1) H_v^{2} \leq B_v^2$ (with $H_v (x)$ denoting the mean curvature at $x$ of the hypersurface $\{y\in\Omega:|v(y)|=|v(x)|\}$), and
$$
|\nabla u_\lambda|^{p-2} |\nabla_{T,u_\lambda} |\nabla u_\lambda||^{2} = \frac{4}{p^2} |\nabla_{T,u_\lambda} |\nabla u_\lambda|^{\frac{p}{2}}|^{2},
$$
we obtain the key inequality 
\begin{equation}\label{StZu2}
\int_{\Omega}\left( \frac{4}{p^2}|\nabla_{T,u_\lambda} |\nabla u_\lambda|^{p/2}|^{2}
+ \frac{n-1}{p-1}H_{u_\lambda}^2 |\nabla u_\lambda|^{p} \right) \eta^2 \, dx
\leq \int_{\Omega}  |\nabla u_\lambda|^{p} |\nabla \eta|^2 \, dx
\end{equation}
for any Lipschitz continuous function $\eta$ with compact support. By taking $\eta=T_s u_\lambda=\min\{s,u_\lambda\}$ in the semistability condition
\eqref{StZu2} we obtain
$$
\int_{\{u_\lambda>s\}}\left( \frac{4}{p^2}|\nabla_{T,u_\lambda} |\nabla u_\lambda|^{p/2}|^{2}
+ \frac{n-1}{p-1}H_{u_\lambda}^2 |\nabla u_\lambda|^{p} \right)\, dx
\leq \frac{1}{s^2}\int_{\{u_\lambda<s\}}  |\nabla u_\lambda|^{p+2}\, dx
$$
for a.e. $s>0$. In particular,
$$
\min\left(\frac{4}{(n-1)p},1\right) I_{p}(u_\lambda-s;\{x\in\Omega:u_\lambda>s\})^p
\leq
\frac{p-1}{(n-1)s^2}\int_{\{u_\lambda<s\}}  |\nabla u_\lambda|^{p+2}\, dx
$$
for a.e. $s>0$, where $I_p$ is the functional defined as follows
\begin{equation*}
I_{p}(v;\Omega):=\left( \int_{\Omega}
\Big(\frac{1}{p'}|\nabla_{T,v} |\nabla v|^{p/2}|\Big)^{2}
+ |H_v|^2 |\nabla v|^p \, dx \right)^{1/p},\quad p \geq 1
\end{equation*}
Then, the $L^\infty$ and $L^r$ estimates of Lemma \ref{lemma1} follow directly from the Morrey and Sobolev type inequalities involving $I_p$ proved in \cite{CaSa}, namely taking $v:=u_\lambda-s$ in
\begin{equation*}
\|v\|_{L^\infty(\Omega)}
\leq C_1|\Omega|^{\frac{p+2-n}{np}}I_{p}(v;\Omega)
\end{equation*}
if $n < p+2$, for some constant $C_1 = C_1 (n,p)$, or
\begin{equation*}
\|v\|_{L^r(\Omega)}
\leq C_2|\Omega|^{\frac{1}{r}-\frac{n-(p+2)}{np}}
I_{p}(v;\Omega)\quad \textrm{for every }1\leq r \leq \frac{np}{n-(p+2)},
\end{equation*}
if $n>p+2$, where $C_2 = C_2 (n,p,r)$. The borderline case $n=p+2$ is slightly more involved, but we are still able to prove a Morrey type inequality as for the case $n <p+2$ (see page 22 in \cite{CaSa} for the details). Lemma \ref{lemma1} is thus proved. \qed

{\bf Sketch of the proof of Lemma \ref{lemma2}.}

Define $\psi(t):= (f(t)-f(0))^{\frac{1}{p-1}}$. By our assumptions we have that $\psi$ is increasing, convex for $t$ sufficiently large and superlinear at infinity. 
Choosing $\varphi=\psi(u_\lambda)$ in the semistability condition \eqref{semistab} and observing that $(p-1)\psi(u_\lambda)^p \psi'(u_\lambda) = f'(u_\lambda) \psi(u_\lambda)^2$, we have
\begin{equation}\label{ned1}
\lambda \int_{\Omega} \psi(u_\lambda)^p \psi' (u_\lambda) \leq \int_{\Omega} |\nabla u_\lambda|^p \psi' (u_\lambda)^2
\end{equation}
On the other hand, multiplying \eqref{problem} by $g(u_\lambda):= \int_{0}^{u_\lambda} \psi'(s)^2 \, ds$ and integrating by parts in $\Omega$ we get:
\begin{equation}\label{ned2}
\int_{\Omega} |\nabla u_\lambda|^p \psi' (u_\lambda)^2 = \lambda \int_{\Omega} (f(u_\lambda) - f(0)) g(u_\lambda) + \lambda f(0) \int_{\Omega}  g(u_\lambda)
\end{equation}
Comparing \eqref{ned1} and \eqref{ned2} it is then easy to see that:
\begin{equation}\label{ned3}
\int_{\Omega} \psi(u_\lambda)^{p-1} h(u_\lambda) \leq f(0) \int_{\Omega}  g(u_\lambda)
\end{equation}
where $h(t) := \int_{0}^{t} (\psi'(t)-\psi'(s)) \psi'(s) \, ds$. Now, thanks to the fact that $h(t) \gg \psi'(t)$ for $t$ large (see page 765 in \cite{BH}), we deduce from \eqref{ned3} that there exists a constant $C$ independent of $\lambda$ such that 
\begin{equation}\label{ned4}
\int_{\Omega} \psi(u_\lambda)^{p-1} \psi' (u_\lambda) \leq C
\end{equation}
In particular, since the asymptotic convexity implies $2 t \psi'(t) \geq \psi (t)$ for any $t$ sufficiently large, from \eqref{ned4} we arrive at
\begin{equation}\label{ned5}
\int_{\Omega} \frac{\psi(u_\lambda)^{p}}{u_\lambda} \leq 2 C
\end{equation}
The desired estimate \eqref{L-p'} is then just a direct consequence of \eqref{ned5} and $p'=p/p-1>1$ since:
$$\int_{\{ u_\lambda > 1 \}} \frac{(f(u_\lambda))^{p'}}{u_\lambda} \, dx  = \int_{\{ u_\lambda > 1 \}} \frac{(f(u_\lambda)-f(0) + f(0))^{p'}}{u_\lambda} \, dx$$
$$\leq \int_{\{ u_\lambda > 1 \}} \frac{(f(u_\lambda)-f(0))^{p'}}{u_\lambda} \, dx + \int_{\{ u_\lambda > 1 \}} \frac{f(0)^{p'}}{u_\lambda} \, dx$$
$$\leq \int_{\Omega} \frac{\psi(u_\lambda)^{p}}{u_\lambda} \, dx + f(0)^{p'} \int_{\{ u_\lambda > 1 \}} \frac{1}{u_\lambda} \, dx \leq 2 C + f(0)^{p'} |\Omega|:= M$$
which proves Lemma \ref{lemma2}.\qed

\end{document}